\documentclass[11pt]{article}
\usepackage{geometry}                
\usepackage{graphicx}
\usepackage{amssymb}
\usepackage{epstopdf}
\usepackage{amsmath,amsthm,amscd,amssymb}
\usepackage{latexsym}
\numberwithin{equation}{section}

\theoremstyle{plain}
\newtheorem{theorem}{Theorem}[section]

\newtheorem{corollary}[theorem]{Corollary}

\theoremstyle{definition}

\theoremstyle{remark}

\newtheorem{case[theorem]}{Case}

\title{Pinned distance sets, Wolff's exponent in finite fields and improved sum-product estimates}
\author{Derrick Hart and Alex Iosevich}

\begin{document}
\maketitle

\begin{abstract} An analog of the Falconer distance problem in vector spaces over finite fields asks for the threshold $\alpha>0$ such that $|\Delta(E)| \gtrsim q$ whenever $|E| \gtrsim q^{\alpha}$, where $E \subset {\Bbb F}_q^d$, the $d$-dimensional vector space over a finite field with $q$ elements (not necessarily prime). Here $\Delta(E)=\{{(x_1-y_1)}^2+\dots+{(x_d-y_d)}^2: x,y \in E\}$. The second listed author and Misha Rudnev (\cite{IR07}) established the threshold $\frac{d+1}{2}$, and in \cite{HIKR07} the authors of this paper, Doowon Koh and Misha Rudnev proved that this exponent is sharp in even dimensions. In this paper we improve the threshold to $\frac{d^2}{2d-1}$ under the additional assumption that $E$ has product structure. In particular, we obtain the exponent $\frac{4}{3}$, consistent with the corresponding exponent in Euclidean space obtained by Wolff (\cite{W99}). 
\end{abstract}           

\section{Introduction} 

\vskip.125in 

\subsection{Distance sets} The classical Erd\H os distance problem asks for the minimal number of distinct distances determined by a finite point set in ${\Bbb R}^d$, $d \ge 2$. The continuous analog of this problem, called the Falconer distance problem asks for the optimal threshold such that the set of distances determined by a subset of ${\Bbb R}^d$, $d \ge 2$, of larger dimension has positive Lebesgue measure. It is conjectured that a set of $N$ points in ${\Bbb R}^d$, $d \ge 2$, determined $ \gtrapprox N^{\frac{2}{d}}$ distances and, similarly, that a subset of ${\Bbb R}^d$, $d \ge 2$, of Hausdorff dimension greater than $\frac{d}{2}$ determines a set of distance of positive Lebesgue measure. Neither problem is close to being completely solved. See \cite{KT04} and \cite{SV05}, and the references contained therein, on the latest developments on the Erd\H os distance problem. See \cite{Erd05} and the references contained therein for the best known exponents for the Falconer distance problem. 

In vector spaces over finite fields, one may define for $E \subset {\Bbb F}_q^d$, 
$$ \Delta(E)=\{||x-y||: x,y \in E\},$$ where 
$$ ||x-y||={(x_1-y_1)}^2+\dots+{(x_d-y_d)}^2,$$ and one may again ask for the smallest possible size of $\Delta(E)$ in terms of the size of $E$. There are several issues to contend with here. First, $E$ may be the whole vector space, which would result in the rather small size for the distance set: 
$$ |\Delta(E)|={|E|}^{\frac{1}{d}}.$$ 

Another annoying consideration is that if $q$ is a prime congruent to $1 \mod(4)$, then there exists $i \in {\Bbb F}_q$ such that $i^2=-1$. This allows us to construct a set 
$$ Z=\{(t,it): t \in {\Bbb F}_q\}$$ and one can readily check that 
$$ \Delta(Z)=\{0\}.$$  

The first result in this direction is proved in \cite{BKT04}. The authors get around the first mentioned obstruction by assuming that $|E| \lesssim q^{2-\epsilon}$ for some $\epsilon>0$. They get around the second mentioned obstruction by mandating that $q$ is a prime $\equiv 3 \mod(4)$. As a result they prove that 
$$ |\Delta(E)| \gtrsim {|E|}^{\frac{1}{2}+\delta},$$ where $\delta$ is a function of $\epsilon$. 

In \cite{IR07} the second author along with M. Rudnev went after a distance set result for general fields in arbitrary dimension with explicit exponents. In order to deal with the obstructions outlined above, they reformulated the question in analogy with the Falconer distance problem: how large does $E \subset {\Bbb F}_q^d$, $d \ge 2$, need to be to ensure that $\Delta(E)$ contains a positive proportion of the elements of ${\Bbb F}_q$. They proved that if $|E| \ge 4q^{\frac{d+1}{2}}$, then $\Delta(E)={\Bbb F}_q$. At first, it seemed reasonable that the exponent $\frac{d+1}{2}$ may be improvable, in line with the Falconer distance conjecture described above. However, the authors along with D. Koh and M. Rudnev discovered in \cite{HIKR07} that the arithmetic of the problem makes the exponent $\frac{d+1}{2}$ best possible in {\bf odd dimensions}, at least in general fields. In even dimensions it is still possible that the correct exponent is $\frac{d}{2}$, in analogy with the Euclidean case. 

The example that shows that the $\frac{d+1}{2}$ is sharp in odd dimensions is very radial in nature and this led the authors of this paper to consider classes of sets that possess a certain amount of product structure. In particular, we shall that if $|E| \subset {\Bbb F}_q^2$ satisfies $|E| \ge Cq^{\frac{4}{3}}$ and $E$ is a product set, then $|\Delta(E)| \ge cq$. This is in line with Wolff's result for the Falconer conjecture in the plane which says that the Lebesgue measure of the set of distances determined by a subset of the plane of Hausdorff dimension greater than $\frac{4}{3}$. In higher dimensions we shall obtain a positive proportion of the distances for products sets of size $\gtrsim q^{\frac{d^2}{2d-1}}$, improving an analog of Erdogan's (\cite{Erd05}) exponent in Euclidean space for general sets. 

\subsection{Pinned distance sets} Let $\pi_j(x)=(x_1, \dots, x_{j-1}, x_{j+1}, \dots, x_d)$ and define 
$$ E^j_z=\pi_j(E) \times \{z\},$$ where $z$ is an element of 
$$ \{z \in {\Bbb F}_q: (x_1,x_2, \dots, x_{j-1}, z, x_{j+1}, \dots, x_d) \in E\}.$$

Define 
$$ \Delta^j_z(E)=\{||x-y||:x \in E, y \in E^j_z\}.$$ 

Our first main result is the following. 
\begin{theorem} \label{distpinned} Let $E \subset {\Bbb F}_q^d$ and let $E^j_z$ be defined with respect to a projection $\pi_j$, for some $1 \leq j \leq d$ and an element $z \in \{z \in {\Bbb F}_q: (x_1,x_2, \dots, z, x_{j+1}, \dots, x_d) \in E\}$ as above. Suppose that 
$$ |E||E^j_z| \ge Cq^d.$$ 

Then 
\begin{equation} \label{orgasm} |\Delta^j_z(E)| \ge q\frac{3C}{3C+1}.\end{equation}
\end{theorem} 

Observe that if $E$ is a product set, then $\Delta^j_z(E) \subset \Delta(E)$. This leads us to the following consequence of Theorem \ref{distpinned}. 
\begin{corollary} \label{productscore} Suppose that $E=E_1 \times E_2 \times \dots \times E_d$, where $E_j$ is contained in ${\Bbb F}_q$. Suppose that 
$$ |E| \ge Cq^{\frac{d^2}{2d-1}}.$$ 

Then 
$$ |\Delta(E)| \ge q\frac{3C}{3C+1}.$$ 
\end{corollary} 

The Corollary follows immediately from Theorem \ref{distpinned} since, after perhaps relabeling some coordinates, we may assume, using straightforward pigeon-holing, that $E=E' \times E''$, where $E' \subset {\Bbb F}_q^{d-1}$ and $|E'| \ge {|E|}^{\frac{d-1}{d}}$. Observe that we could have made a much weaker, though more technical, assumption on the structure of $E$. We shall attempt to classify this notion in a precise way in a subsequent paper. 

\vskip.125in 

\subsection{Sums and products} A related question that has recently been attacked using a similar Fourier-geometric framework is the problem of sums and products in ${\Bbb F}_q$ in the following form. Let $A \subset {\Bbb F}_q$. How large does $A$ need to be to ensure that 
$$ {\Bbb F}_q^{*} \subset A \cdot A+A \cdot A,$$ or, more modestly, 
$$ |A \cdot A+A \cdot A| \ge cq$$ for some $c>0$. 

The authors of this paper have recently proved the following result. 
\begin{theorem} \label{vzhopu} Let $A \subset {\Bbb F}_q$. \begin{itemize} \item If $|A|>q^{\frac{3}{4}}$,
then 
$$ {\Bbb F}_q^{*} \subset A \cdot A+A \cdot A.$$ 

\item If $|A| \ge C^{\frac{1}{2}}_{size}q^{\frac{2}{3}}$, then 
$$ |A \cdot A+A \cdot A| \ge q \frac{C^{\frac{3}{2}}_{size}}{1+C^{\frac{3}{2}}_{size}}.$$ 
\end{itemize} 
\end{theorem} 

We use the method of proof of Theorem \ref{distpinned} above to obtain the following. 
\begin{theorem} \label{pinnedsp} Let $E\subset \mathbb F_q^d$ and let $E^j_z$ be defined as above. Then
\begin{equation} \label{orgasm2} |\{x\cdot y: x\in E_z, y \in E\}| \ge q \frac{2C}{2C+1} \end{equation} whenever 
$$ |E^j_z||E| \ge Cq^d.$$ 
\end{theorem}

Setting $E=A \times A \times A$, we obtain the following consequence recently obtained by Shparlinski (\cite{S07}) in the case $d=2$. 

\begin{corollary} Let $A,B$ be subsets of ${\Bbb F}_q$. Suppose that 
$$ |A| \ge Cq^{\frac{2}{3}}.$$ 

Then for any $z \in A$, 
$$ |A \cdot A+zA| \ge \frac{2C}{2C+1}q.$$ 
\end{corollary} 

\vskip.125in 

\section{Proof of Theorem \ref{distpinned}}

\vskip.125in 

We drop the exponent $j$ in the sequel for the sake of convenience. Define 
$$ \nu(t)=|\{(x,y) \in E_z \times E: ||x-y||=t\}|,$$ where, as usual 
$$ ||x||=x_1^2+\dots+x_d^2.$$ 

By Cauchy-Schwartz, 
$$ \nu^2(t) \leq |E_z| \cdot \sum_{||x-y||=||x-y'||=t} E_z(x)E(y)E(y'),$$ so 
$$ \sum_t \nu^2(t) \leq |E_z| \cdot \sum_{||x-y||=||x-y'||} E_z(x)E(y)E(y')$$
$$=q^{-1}|E_z| \cdot \sum_s \sum_{x,y,y'} \chi(s(||x-y||-||x-y'||)) E_z(x)E(y)E(y')$$
$$=q^{-1}{|E_z|}^2{|E|}^2+R,$$ and 
$$ R=q^{-1}|E_z| \sum_{s \not=0} \sum_{x \in E_z} 
{\left| \sum_{y \in E} \chi(s(||y||-2x \cdot y)) \right|}^2,$$ since 
$$ ||x-y||-||x-y'||=(||y||-2x \cdot y)-(||y'||-2x \cdot y').$$ 

It follows that 
$$ R \leq q^{-1}|E_z| \sum_{s \not=0} \sum_{x \in {\Bbb F}_q^{d-1} \times \{z\}} \sum_{y,y' \in E} 
\chi(-2sx \cdot (y-y')) \chi(s(||y||-||y'||))$$
$$=q^{d-2}|E_z| \sum_{s \not=0} \sum_{\pi_d(y)=\pi_d(y')} E(y)E(y') \chi(-2sz(y_d-y'_d)) 
\chi(s(y_d^2-{y'}_d^2))$$
$$=q^{d-2}|E_z| \sum_{s} \sum_{\pi_d(y)=\pi_d(y')} E(y)E(y') \chi(-2sz(y_d-y'_d)) 
\chi(s(y_d^2-{y'}_d^2))$$
$$-q^{d-2}|E_z| \sum_{\pi_d(y)=\pi_d(y')} E(y)E(y')=A-B.$$

Now, 
$$ B \leq q^{d-2}|E_z||E|q=q^{d-1}|E_z||E|$$ because one coordinate can contribute at most a factor of $q$. On the other hand, 
$$ A=q^{d-1}|E_z| \cdot \sum_{2z(y_d-y'_d)=y_d^2-{y'}_d^2; \pi_d(y)=\pi_d(y')} E(y)E(y')$$
$$=q^{d-1}|E_z| \cdot \sum_{2z=y_d+y'_d; y_d \not=y'_d; \pi_d(y)=\pi_d(y')} E(y)E(y')$$
$$+q^{d-1}|E_z| \sum_{y=y'} E(y)E(y')$$
$$ \leq 2q^{d-1}|E_z||E|.$$

\subsection{Conclusion of the proof:} We have 
$$ {|E|}^2{|E_z|}^2={\left( \sum_t \nu(t) \right)}^2 \leq |\Delta_z(E)| \cdot \sum_t \nu^2(t)$$
$$ \leq |\Delta_z(E)| \cdot ({|E|}^2{|E_z|}^2q^{-1}+3q^{d-1}|E_z||E|).$$ 

It follows that 
$$ |\Delta_z(E)| \ge cq$$ if 
$$ |E||E_z| \ge Cq^d.$$ 

\vskip.25in 

\section{Proof of Theorem \ref{pinnedsp}}

\vskip.125in 

Consider the incidence function 
$$\nu(t)=|\{(x,y)\in E_z \times E: x\cdot y =t\}|.$$
Then
$$\sum_t \nu^2(t) \leq |E_z|\sum_{x\cdot y=x\cdot y'} E_z(x)E(y)E(y')$$
$$=|E_z|^2|E|^2q^{-1}+|E_z|q^{-1}
\sum_{s\neq 0}\sum_{\substack{x \in E_z \\ y,y' \in E}} \chi(sx\cdot(y-y'))$$
$$=I+II$$
Now
$$II = |E_z|q^{-1}\sum_{s\neq 0}\sum_{x \in E_z}\left| \sum_{y\in E} \chi(sx\cdot y)\right|^2$$
$$\leq |E_z|q^{-1}\sum_{s\neq 0}
\sum_{x \in \mathbb F_q^{d-1}\times \{z\} }\left| \sum_{y\in E} \chi(sx\cdot y)\right|^2$$
$$= |E_z|q^{d-2}\sum_{s\neq 0}
\sum_{\pi_d(y)=\pi_d(y')} \chi(sz(y_d-y_d'))$$
$$= |E_z||E|q^{d-1}
-|E_p|q^{d-2}\sum_{\pi_d(y)=\pi_d(y')} 1$$
$$=A-B.$$

We conclude that 
$$B \leq  |E_z||E|q^{d-1},$$ which implies that 
$$|II|\leq 2 |E||E_z|q^{d-1}$$
and 
$$\sum_t \nu^2(t) \leq |E_z|^2|E|^2q^{-1}+2|E||E_z|q^{d-1},$$ and the conclusion follows in the same way as in the proof of Theorem \ref{distpinned}.

\enddocument